\documentclass[10pt]{article}
\usepackage{amsmath,amssymb}
\newcommand{\dis}{\displaystyle}
\newcommand{\pf}{\noindent{\bf Proof. }}
\newcommand{\dint}{\int_{\R^3}}

\newcommand{\Div}{\mathrm{div}}
\newcommand{\dt}{{\partial_t}}

\newcommand{\Mp}{{\mathcal M}^+(\R^3)}
\newcommand{\R}{\mathbb{R}}

\newcommand{\un}{{\mathbf{1}}}
\newcommand{\into}{\int_\Omega}
\newcommand{\unk}{{\mathbf{1}}_{\{|u|\geq C_k\}}}
\def\na{\nabla}
\def\pa{\partial}

\def\qed{\hbox{${\vcenter{\vbox{
  \hrule height 0.4pt\hbox{\vrule width 0.4pt height 6pt
  \kern5pt\vrule width 0.4pt}\hrule height 0.4pt}}}$}}

\newtheorem{theo}{Theorem}
\newtheorem{prop}[theo]{Proposition}
\newtheorem{lemm}[theo]{Lemma}
\newtheorem{coro}[theo]{Corollary}

\author{Antoine Mellet\thanks{University of British Columbia, Department of Mathematics, 1984 Mathematics Road, Vancouver, BC}, Alexis Vasseur \thanks{University of Texas at Austin, department of Mathematics,
1 University Station C1200, Austin, TX 78712}}
\title{$L^p$ estimates for quantities advected by a compressible flow}
\date{}

\begin{document}
\maketitle

\bibliographystyle{plain}

 {\small \noindent{\bf Abstract:}
We consider the evolution of a quantity advected by a compressible flow and subject to diffusion.
When this quantity is scalar it can be, for instance, the temperature of the flow or the concentration of some pollutants.
Because of the diffusion term, one expects the equations to have a regularizing effect. However,
in their Euler form, the equations describe the evolution of the quantity multiplied by the density of the flow. The parabolic structure is thus degenerate near vacuum (when the density vanishes).
In this paper we show that we can nevertheless derive uniform $L^{p}$ bounds that do not depend on the density (in particular the bounds do not degenerate near vacuum).
Furthermore the result holds even when the density  is only a measure.

We investigate both the scalar and the system case. In the former case, we obtain $L^{\infty}$ bounds.
In the latter case the quantity being investigated could be the velocity field in compressible Navier-Stokes type of equations, and we derive uniform $L^p$ bounds for some $p$ depending on the ratio between the two viscosity coefficients (the main additional difficulty in that case being to deal with the second viscosity term involving the divergence of the velocity).
Such estimates are, to our knowledge, new and interesting since they are uniform with respect to the density.
The proof relies mostly on a method introduced by De Giorgi to obtain regularity results for elliptic equations with discontinuous diffusion coefficients. }

\section{Introduction}
Let $\theta(t,x)$ be a function defined on
$[0,T]\times\R^3$, solution to the following equation:
\begin{equation}\label{eq:scalaire}
\begin{array}{l}
\dt (\rho\,  \theta) +\Div(\rho\, v\,\theta)-\Div (\mu \nabla \theta)=\rho
F+\Div (\rho\, G),\\
\theta(0,x)=\theta_0(x),\\
\end{array}
\end{equation}
where $F$, $G_1$, $G_2$, $G_3$, $\rho$, and $v$ are given functions such that
$(\rho, v)$ satisfies the following continuity equation:
\begin{equation}\label{eq:transport_rho}
\begin{array}{l}
\dt \rho+\Div (\rho\, v)=0,\\
\rho(0,x)=\rho_0(x).
\end{array}
\end{equation}
At this point, we wish to stress out the fact that we will be considering very general diffusion
coefficients $\mu(t,x)$ throughout the paper: We will only assume that $\mu(t,x)$ is   measurable and verifies
\begin{equation}\label{eq:mu}
\mu(t,x)\geq 1 \qquad \mathrm{for} \ \ t\in [0,T], \ x\in \R^3.
\end{equation}

Such a system of equations arises in a lot of contexts.
The function $\theta$ can, for instance, model the temperature of a fluid with density $\rho$
and velocity $v$, or the density of pollutant spreading in this
fluid.
When $\mu$ satisfies (\ref{eq:mu}),  it is well known that the usual energy inequality gives some bounds on $\theta$ (typically in $L^2(0,T,L^6(\R^3))$).
However no further $L^p$ estimates can be obtained on $\theta$ directly. This is because the conserved quantities always involve the density $\rho$ (they are typically of the form $\rho\theta^\alpha$),
so that any bounds obtained by this mean will degenerate on cavities (i.e. when $\rho=0$).
The goal of this paper is to establish $L^\infty$ estimates on $\theta$ that do
not degenerate on cavities and that do not depend too much on the density $\rho$.

Note that this kind of estimates is quite natural if one looks at the evolution of $\theta$ in the Lagrangian description of the flow.
However, working in the Lagrangian framework would be of little help here, since all the bounds on $F$ and $G$ would be tampered (and we would lose the divergence structure of the term $\Div(\rho\, G)$).

\vspace{10pt}

Let us now state our results more precisely.
First, we  will be  considering only "suitable" solutions of (\ref{eq:scalaire}).
By this, we mean functions that verify the following inequality:
\begin{equation}\label{eq:energiescalaire}
\frac{d}{dt}\int_{\R^3}\rho \phi(\theta)\,dx+\int_{\R^3}\mu
\phi''(\theta)|\nabla \theta|^2\,dx\leq \int_{\R^3}\rho
F\phi'(\theta)\,dx+\dint\phi'(\theta)\Div (\rho G)\,dx,
\end{equation}
for every convex function $\phi\in W^{2,\infty}_{\mathrm{loc}}(\R)$ verifying:
\begin{equation}\label{eq:phi}
\lim_{y\to\infty} \frac{\phi(y)}{y^2}\leq 1.
\end{equation}
As we will see in Section 2, any regular solution of (\ref{eq:scalaire}), (\ref{eq:transport_rho}) verifies
(\ref{eq:energiescalaire}) (with an equality).
It is thus very natural to consider weak solutions verifying (\ref{eq:energiescalaire}), in the same spirit as that of Leray's weak solutions for incompressible Navier-Stokes equations.

We denote by $\Mp$ the set of positive measures in $\R^3$, and for any
$\rho\in \Mp$, we denote  by $L^2(\rho)$ the set of $\rho$-measurable
function $h$ satisfying $\int h^2 d\rho<\infty$.
The first result of this paper is the following:
\begin{theo}\label{theo:scalaire}
Take $T$ finite or $T=+\infty$.
Assume that the viscosity coefficient $\mu$ verifies
(\ref{eq:mu}), that $\rho$ lies in  $L^\infty(0,T;\Mp)$
 and that $F$ and $G$ are such that there
exists $0<\alpha<1$ such that
\begin{equation}\label{eq:FGLinfty}
\rho^\alpha|F|+\sum_{i=1}^3\rho^{1+\alpha}|G_i|^2/\mu\in
L^p(0,T;L^q(\R^3)),
\end{equation}
for some $p$ and $q$ satisfying:
\begin{eqnarray*}
&&p>\frac{1}{1-\alpha},\qquad \frac{2}{p}+\frac{3}{q}< 2.
\end{eqnarray*}
Let $\theta\in L^\infty(0,T;L^2(\rho(t)))$ with $\nabla\theta\in
L^2((0,T)\times\R^3))$ be a solution of (\ref{eq:energiescalaire}) for
every $\phi\in W^{2,\infty}_{\mathrm{loc}}(\R)$ verifying (\ref{eq:phi}).

Then the following results hold:
\begin{itemize}
\item If $\theta_0\in L^\infty(\R^3)$, then $\theta\in
L^\infty([0,T]\times\R^3)$.

\item If $\theta_0$ is only bounded in $L^2(\rho_0)$ and if $\rho\in L^\infty(0,T;L^r(\R^3))$
for some $r>3/2$, then $\theta\in L^\infty(t_0,T;L^\infty(\R^3))$
for every $t_0>0$.
\end{itemize}
\end{theo}
\vspace{10pt}

One of the main motivation of this article is to derive some bounds for the velocity field of compressible flows (see the motivations subsection below).
In the compressible Navier-Stokes system of equations, the velocity is advected by itself and subject to  viscosity effects, as previously.
However, there are now two viscosity terms, one of which involves $\Div\, u$ which induces a strong system structure.
Our first result does not apply to that case, so we will also study the system case:
Consider a vector-valued function
$u\in\R^3$ solution in $[0,T]\times\R^3$ to the following system of equations:
\begin{eqnarray}
&&\dt \rho u +\Div(\rho\,  v\otimes u)-\Div (2\mu \nabla
u)-\nabla(\lambda\Div u)=\rho F+\Div(\rho\,
G),\label{eq_NS2}\\
&&u(0,x)=u_0(x),
\end{eqnarray}
where $(\rho,v)$ still verifies (\ref{eq:transport_rho}).
In addition to (\ref{eq:mu}), we will assume that the second viscosity
coefficient $\lambda$ (which can also depends on $(t,x)$) verifies
for every $(t,x)\in[0,T]\times\R^3$:
\begin{eqnarray}
&&\nu(t,x)=2\mu(t,x)+3\lambda(t,x)\geq1, \label{eq:h+Ng}\\
&&3|\lambda(t,x)|\leq \kappa \nu(t,x),\label{eq:lambda}
\end{eqnarray}
for some $0<\kappa<1/2$.

As in the scalar case, we  will consider
"suitable" solutions of (\ref{eq_NS2}), which verify:
\begin{equation}\label{eq_NSenergy}
\begin{array}{l}
\qquad\dis{ \frac{d}{dt}\int_{\R^3}\rho
\phi(|u|)\,dx}\\[0.3cm]
\dis{+\int_{\R^3} \nu\frac{\phi'(|u|)}{|u|}|\nabla u|^2\,dx+
\int_{\R^3} \nu\left[\phi''(|u|)-\frac{\phi'(|u|)}{|u|}\right]|\nabla |u||^2\,dx}\\[0.3cm]
\dis{\leq -\int_{\R^3}\lambda
\left[\phi''(|u|)-\frac{\phi'(|u|)}{|u|}\right]\sum_{ij}\pa_i
u_j\frac{u_iu_j}{|u|^2}(\Div u)\,dx+ \int_{\R^3}\rho\frac{u}{|u|}\cdot F\phi'(|u|)\,dx}\\[0.3cm]
\dis{+ \dint\phi'(|u|)\sum_{i=1}^3\sum_{j=1}^3
\frac{u_j}{|u|}\partial_i(\rho G_{ij})\,dx},
\end{array}
\end{equation}
for every convex function $\phi\in
W^{2,\infty}_{\mathrm{loc}}(0,\infty)$ verifying:
\begin{equation}\label{eq:phi2}
\phi''(y)-\phi'(y)/y\geq0 \qquad\lim_{y\to\infty}
\frac{\phi(y)}{y^2}\leq 1.
\end{equation}
Again, we will see in section 2  that any regular solution of
(\ref{eq_NS2}), (\ref{eq:transport_rho}) verifies (\ref{eq_NSenergy}).

The system structure weakens the bounds, and our  second result is the following:
\begin{theo}\label{theo:systeme}
Take $T$ finite or $T=+\infty$.
Assume that $\rho\in L^\infty(0,T;\Mp)$, that  $\mu$ and $\lambda$
verify (\ref{eq:mu}), (\ref{eq:h+Ng}) and (\ref{eq:lambda}), and that $F$ and $G$
 are such that
\begin{equation}\label{eq:FGLinfty2}
\sum_{i}\rho^\alpha|F_i|+\sum_{i,j}\rho^{1+\alpha}|G_{ij}|^2/\nu\in
L^p(0,T;L^q(\R^3)),
\end{equation}
for some $p$ and $q$ satisfying
\begin{eqnarray}\label{eq:pq}
&&p>\frac{1}{1-\alpha}, \qquad \frac{2}{p}+\frac{3}{q}< 2.
\end{eqnarray}
Consider $u\in
L^\infty(0,T;L^2(\rho(t)))$, with $\nabla u\in
L^2((0,T)\times\R^3)$, solution to (\ref{eq_NSenergy}) for every
$\phi\in W^{2,\infty}_{\mathrm{loc}}(0,\infty)$ verifying
(\ref{eq:phi2}). Then the two following results hold true:
\begin{itemize}
\item If $u_0\in L^\infty(\R^3)$, then for any ball $B\subset \R^3$, $u\in
L^{2+\log_2(1/\kappa)}([0,T]\times B)$ where $\kappa$ is defined in
(\ref{eq:lambda}).

\item If $u_0$ is only bounded in $L^2(\rho_0)$ and $\rho\in L^\infty(0,T;L^r(\R^3))$
for some $r>3/2$, then $u\in
L^{2+\log_2(1/\kappa)}_{\mathrm{loc}}((0,T]\times\R^3)$.
\end{itemize}
\end{theo}
Note that the proofs of both theorems will only make use of the relations
(\ref{eq:energiescalaire}) or (\ref{eq_NSenergy}) which do not
depend on the advection velocity $v$. In particular the results hold
true even if we cannot give a meaning to (\ref{eq:scalaire}),
(\ref{eq:transport_rho}) or (\ref{eq_NS2}). This explains why no
assumption is required on $v$. Let us aslo emphasize that no lower
bound is required on the density $\rho$, and that we can also deal
with a density which is merely a measure.
Note also that when $\rho\equiv1$, $v=0$, $\mu\equiv1$ and $\lambda=0$, the equation (\ref{eq:scalaire}) is
nothing but the heat equation with a given right hand side.
In that case,  the conditions on $p$ and $q$ in (\ref{eq:FGLinfty}) and
(\ref{eq:FGLinfty2}) are the usual one to obtain $L^\infty$
regularity of the solutions of parabolic equations.
Our conditions are thus  optimal in that sense.
Finally, we stress out the fact that these conditions also allow us to consider forces
that blow up near vacuum like $1/\rho^\alpha$ for some $0\leq\alpha<1$.
\vspace{10pt}

In the next subsections, we show how those results apply to compressible
Navier-Stokes equations and we give the main idea of the proof.
It relies on a method introduced by De Giorgi in \cite{DeGiorgi} to
show $C^\alpha$ regularity of solutions to elliptic equations with
rough diffusion coefficients (measurable).
This method was used for the first time on Navier-Stokes equation in
\cite{Vasseur}, where  an alternative proof of partial regularity for
solutions to incompressible Navier-Stokes equation, first proven by
Caffarelli, Kohn and Nirenberg in \cite{CKN}, was provided.
Note that the method makes use  of inequalities (\ref{eq:energiescalaire}) and
(\ref{eq_NSenergy}) which describe the evolution of quantities of the
form:
$$
\dint \rho \phi(\theta)\,dx, \qquad \dint \rho \phi(|u|)\,dx.
$$
This idea was already the key stone of the paper \cite{MV}.

\subsection{Motivation: Compressible Navier-Stokes equation}

In this subsection we describe the consequences of our results on compressible barotropic Navier-Stokes equations.
Our aim is to show the pertinence of the result but also its limits.
We consider the following system of equations:
\begin{equation}\label{eq_NSvrai}
\begin{array}{l}
 \pa_{t} \rho +\Div (\rho u) = 0 \\
 \pa_{t}(\rho u) +\Div (\rho u \otimes u) + \nabla_{x} \rho^\gamma -\Div (2\mu\, \na u) - \na (\lambda\, \Div u)=
 0,
\end{array}
\end{equation}
where the viscosity coefficients verify:
\begin{equation}
 \mu\geq1 \qquad 2\mu+3\lambda \geq1.
\end{equation}
Note that this system can be written in  the form of (\ref{eq_NS2}), (\ref{eq:transport_rho}) with $v=u$, $F=0$ and $G=-\rho^{\gamma-1}$.

The problem of the existence of solutions defined globally in time for this type of
system was addressed in one dimension for smooth enough data by
Kazhikov and Shelukhin \cite{K}, and for discontinuous one, but
still with densities away from zero, by Serre \cite{Serre} and Hoff
\cite{Hoff}. Those results have been generalized to higher
dimensions by Matsumura and Nishida \cite{nishida} for smooth data
close to equilibrium and by Hoff \cite{Hoff2}, \cite{Hoff3} in the
case of discontinuous data.

Concerning large initial data, Lions showed in \cite{PLL2} the
global existence of weak solutions for $\gamma\geq3/2$ for $N=2$ and
$\gamma\geq 9/5$ for $N=3$. This result has been extended later by
Feireisl, Novotny, and Petzeltova   to the range $\gamma>3/2$ in
\cite{F2}. 
Other results provide the full range $\gamma>1$ under symmetries
assumptions on the initial datum (see for instance Jiang and Zhang
\cite{J1}). Notice that all those results hold with constant
viscosity coefficients $\mu$ and $\lambda$. Unfortunately, those
theories do not provide enough $L^p$ bounds on the pressure
$\rho^\gamma$ to apply Theorem \ref{theo:systeme}. Indeed, Vaigant
in \cite{Vaigant} showed that, even with rather smooth data, we
cannot expect such bound for $p$ large. Hence Theorem
\ref{theo:systeme} provides only the following partial result:
\begin{coro}
Let $\mu$ and $\lambda$ be two constants verifying (\ref{eq:h+Ng}) and
(\ref{eq:lambda}), and let $\gamma>3/2$. Let $(\rho,u)$ be the solution
of (\ref{eq_NSvrai}) constructed in \cite{F2}. If $T$ is such that $\rho\in L^\infty(0,T;L^p(\R^3))$ with  $p>3\gamma$, then
$$
u\in L^{2+\log_2(1/\kappa)}_{\mathrm{loc}}((0,T]\times\R^3).
$$
\end{coro}
Note that the condition $\gamma>3/2$ is necessary to make use of \cite{F2}. Naturally, the initial data is assumed to have finite energy, which, in particular guarantees that $u_0$ is bounded in $L^2(\rho_0)$.
\vspace{10pt}

All the previously mentioned results only hold for constant viscosity coefficients. From a physical point of view, however, the viscosity coefficients $\lambda$ and $\mu$ are known
to depend on the temperature and thus, for the barotropic model, on
the density.
Unfortunately, very little is known on the existence of
solutions in this case (the main difficulty being to get some compactness on
the density $\rho$).
Bresch and Desjardins have found a new
mathematical entropy for a class of such coefficients which gives
some norms on the gradient of $\rho$ (see \cite{BD3},\cite{BD2}).
This allowed them to construct solutions of (\ref{eq_NSvrai}) when
additional physical terms are added (drag force, Korteweg type term
or cold pressure).
It was shown in \cite{MV} that these additional terms can be droped.
However Theorem \ref{theo:systeme}
cannot be applied in this framework since the new $BD$ entropy
requires to violate condition (\ref{eq:mu}). Thus, in this context, our
results only give a priori bound for a system whose solutions are not yet
known to exist:
\begin{coro}(a priori bounds)
Let $\gamma>3/2$, $\mu(\rho)$ and $\lambda(\rho)$ verify
(\ref{eq:h+Ng}) and (\ref{eq:lambda}) with
$\nu(\rho)\geq\rho^{\beta}$ for some $\beta>4\gamma/3$. Then any
suitable solution $(\rho,u)$ of (\ref{eq_NSvrai}), verifies:
$$
u\in L^{2+\log_2(1/\kappa)}_{\mathrm{loc}}((0,T]\times\R^3).
$$
\end{coro}
Note that in this corollary no assumption needs to be made on
$\rho$. Indeed we take advantage of the fact that the condition on
$G$ depends on $\nu$: We have
\begin{eqnarray*}
\rho^{1+\alpha}\frac{|G|^2}{\nu}&=&\frac{|\rho^{\gamma}|^2}{\rho^{(1-\alpha)}\nu}\\
&\leq &\rho^{2\gamma-1+\alpha-\beta}
\end{eqnarray*}
which is dominated by a power of $\rho$ slightly better than
$2/3\gamma$. So we can apply the theorem, using the fact that the usual entropy
inequality gives $\rho$ bounded in  $L^\infty(L^\gamma)$.

When considering a $\rho$-dependent viscosity coefficient
$\mu(\rho)$, another interesting problem arises. Indeed, the actual form of
the first viscosity term should then be:
$$
\Div (\mu(\rho) D(u)),
$$
where $D(u)=\nabla u+(\nabla u)^T$ is the symmetric part of the
gradient of $u$. Note that when $\mu$ is constant this remark is not relevant
since we  have:
$$
\Div(\mu D(u))=\Div(\mu\na u)+\mu\sum_{i}\partial_{ij}u=\Div(\mu\na
u)+\na(\mu\Div u).
$$
For our purpose, working with $D(u)$ instead of $\na u$
would be far more tedious, since it would add another pure system term.
Finding a equivalent of Theorem \ref{theo:systeme} with such a term
is therefore a very interesting and challenging open problem.

To conclude this remarks,
let us mention that the same type of corollaries  (with the same restrictions) can be applied to the Navier-stokes equation
with temperature.
For instance,
in the recent existence result of
Feireisl \cite{F3}, the viscosity coefficients are
allowed to depend on the temperature and are required to verify the condition
(\ref{eq:mu}).
Let us also mention the result of Bresch and
Desjardin \cite{BDTemp} in which density dependent viscosity
coefficients were considered. However, in this last paper, the
coefficients have to vanish on vacuum, so our results do not apply.

\subsection{Idea of the Proof}
In this subsection we want to describe the main idea of the proof in a simplified framework.
We thus focus on the first statement of Theorem \ref{theo:scalaire} and assume that  $G=0$ and $\rho^{1/3}F$ is
bounded.

We introduce a sequence of  functions $\phi_k$:
$$
\phi_k(y)=[|y|-C_k]^2_+,
$$
where $[z]_+=\sup(0,z)$ and $C_k$ is an increasing sequence of
number defined (in this particular example) by:
$$
C_k=K(1-2^{-k}),
$$
where $K$ will be chosen later.
(Note that $C_k$ converges to  $K$
when $k\to\infty$). With these notations, we set:
$$
U_k=\sup_{0\leq t\leq T}\left(\int_{\R^3}\rho
\phi_k(\theta)\,dx\right)+\int_0^T\int_{\R^3}\mu
\phi_k''(\theta)|\nabla \theta|^2\,dx\,dt.
$$
If we think of $\int\rho \theta^2\,dx$ as an energy (as in the case of compressible Navier-Stokes equation), then we can think of $\int\rho
\phi_k(\theta)\,dx$ as a level set of the energy, namely the
energy corresponding to the values of $|\theta|$ that are greater than $C_k$.
The quantity $U_k$ is thus the sum of the supremum of the $k$-level set of energy
on all the times $0\leq t\leq T$ and the viscous
dissipation of this $k$-level set of energy over the same time interval.
Our goal is now to determine how
the quantity $U_k$ depends on the previous
$(k-1)$-level set quantity $U_{k-1}$.

This is done by using the inequality (\ref{eq:energiescalaire})
with $\phi=\phi_k$ and integrating it
over $[0,t]$ for every $0\leq t\leq T$.
Noticing that if
$K>2\|\theta_0\|_{L^\infty}$, then
$$
\dint \rho_0\phi_k(\theta_0)\,dx=0 \qquad \mbox{ for all $k>1$},
$$
we deduce:
\begin{eqnarray*}
U_{k}&\leq&
2\int_0^T\dint|\rho F\phi_k'(\theta)|\,dx\,dt\\
&\leq&C\int_0^T\dint\rho^{2/3}|\phi_k'(\theta)|\,dx\,dt.
\end{eqnarray*}
The next step is to control the right-hand side in terms of
$U_{k-1}$. The main ingredients are Sobolev embedding and Tchebychev
inequality. We denote:
$$
\theta_{k}=(|\theta|-C_k)_+.
$$
Since $U_{k-1}$ controls the square of the
$L^\infty(0,T;L^2(\R^3))$ norm of $\rho^{1/2} \theta_{k-1}$ together with the
square of the $L^2([0,T]\times\R^3)$ norm  of the gradient in $x$ of
$\theta_{k-1}$, Sobolev embeddings and H\"older inequalities yield:
\begin{equation}\label{Sobolev}
\|\rho^{1/5}\theta_{k-1}\|^2_{L^{10/3}([0,T]\times\R^3)}\leq
CU_{k-1}.
\end{equation}

Next, using  (\ref{Sobolev}) and the fact that when $|\theta|\geq C_k$, we have $\theta_{k-1}\geq
C_k-C_{k-1}$ and so
$1_{\{\theta_k>0\}}\leq\frac{\theta_{k-1}}{C_k-C_{k-1}}$, we get:
\begin{eqnarray*}
&&\qquad\int_0^T\dint\int_{\R^3}\rho^{2/3}|\phi_k'(\theta)|\,dx\,dt\\
&&=2\int_0^T\dint \rho^{2/3}\theta_k
{\mathbf{1}}_{\{|\theta|\geq C_k\}}\,dx\,dt\\
&&\leq 2\int_0^T\dint \rho^{2/3}\theta_{k-1}
{\mathbf{1}}_{\{\theta_{k-1}\geq C_k-C_{k-1}\}}\,dx\,dt\\
&&\leq 2\int_0^T\dint\rho^{2/3}
\theta_{k-1}\frac{\theta^{7/3}_{k-1}}{(C_k-C_{k-1})^{7/3}}
{\mathbf{1}}_{\{\theta_{k-1}\geq C_k-C_{k-1}\}}\,dx\,dt\\
&&\leq \frac{C}{(C_k-C_{k-1})^{7/3}}U_{k-1}^{5/3},
\end{eqnarray*}
which leads to:
$$
U_k\leq \frac{C2^{7k/3}}{K^{7/3}}U_{k-1}^{5/3}.
$$
Note that in more general cases we will get an estimate of the form:
$$
U_k\leq \frac{C2^{\alpha k}}{K^{\underline{\gamma}}}U_{k-1}^{\beta}.
$$
What is important in this inequality is to have $\beta>1$ and $\underline{\gamma}>0$.
As a matter of fact, we can then prove (see Lemma \ref{lemm_technique_scalar}) that for any given
$U_1$, there exists a $K$ large enough for which any sequence
$U_k$ satisfying this induction inequality will converge to $0$.
The definition of $U_k$ then yields
 that $\na
\theta_\infty$ is zero, with $\theta_\infty=(\theta-K)_+$. This implies that
$\theta_\infty$ is constant and this constant has to be 0 since
$\theta_\infty\in L^2(L^6)$.
It follows that  $\theta\leq K$ almost everywhere.
\vspace{10pt}

Note that the "non-linearisation" process, which gives a $\beta>1$,
is primordial to counter fight the growth of $2^{\alpha k}$. This is
provided by the double action of the (linear) Sobolev embedding
which gives $L^p$-norms with $p>2$ and the non-linear Tchebychev
type inequality which can be used thanks to the fact that we
consider a previous energy level set. This is the  key idea
of De Giorgi's method.

When considering an initial value which is not bounded, we need to
"non-linearize" the contribution of $\int\rho\theta^2\,dx$ near
$t=0$. This can be done introducing increasing time $T_k$ and
integrating (\ref{eq:energiescalaire}) on $[T_k,T]$.

The difficulties considering the system case are more serious. The
problem is to control the viscous term of the right-hand side part
of (\ref{eq_NSenergy}). This can be done only in a linear way. This
forces to take a sequence $C_k$ converging to infinity. The $L^p$
norm will then depends on the rate of decreasing of $U_k$ from a
level set to the next.

\section{Consistency of the notion of suitable solutions}

In this short section we check that the notion of suitable solutions
is consistent with the differential equations.
This is quite important for us since the method
relies entirely on the estimates (\ref{eq:energiescalaire}) and
(\ref{eq_NSenergy}). We will show the two following lemmas:
\begin{lemm}\label{lemm_consistency_scalaire}
Let $\theta$ be a regular solution to (\ref{eq:scalaire})
(\ref{eq:transport_rho}), then $\theta$ satisfies
(\ref{eq:energiescalaire}) with an equality.
\end{lemm}
\begin{lemm}\label{lemm_consistency_NS}
Let $u$ be a regular solution to (\ref{eq_NS2})
(\ref{eq:transport_rho}), then $u$ satisfies (\ref{eq_NSenergy}) for
any $\phi$ verifying (\ref{eq:phi2}).
\end{lemm}
\noindent{\bf Proof of Lemma \ref{lemm_consistency_scalaire}:} Using
(\ref{eq:scalaire}) and (\ref{eq:transport_rho}), we find:
$$
\rho\dt \theta +\rho v\cdot\nabla\theta-\Div (\mu \nabla
\theta)=\rho F+\Div (\rho G).
$$
Multiplying by $\phi'(\theta)$ and integrating in $x$, we find:
\begin{eqnarray*}
&&\dint \rho\dt \phi(\theta)\,dx+\dint \rho
v\cdot\nabla\phi(\theta)\,dx+\dint \mu\phi''(\theta) |\nabla
\theta|^2\,dx\\
&&\qquad=\dint \rho F\phi'(\theta)\,dx+\dint \phi'(\theta)\Div (\rho
G)\,dx.
\end{eqnarray*}
Using again (\ref{eq:transport_rho}) gives the result.\qquad \qed
\vskip1cm \noindent{\bf Proof of Lemma \ref{lemm_consistency_NS}:}
>From (\ref{eq_NS2}) and (\ref{eq:transport_rho}) we find:
$$
\rho\dt u +\rho v\cdot\nabla u-\Div (2\mu \nabla
u)-\nabla(\lambda\Div u)=\rho F+\Div (\rho G).
$$
Multiplying it by $\phi'(|u|)u/|u|$ and integrating with respect to
$x$ we find:
\begin{eqnarray*}
&&\qquad\dint\rho\dt(\phi(|u|))\,dx+\dint \rho v \cdot
\nabla(\phi(|u|))\,dx+\dint 2\mu\frac{\phi'(|u|)}{|u|}|\nabla
u|^2\,dx\\
&&+\dint \lambda\frac{\phi'(|u|)}{|u|}|\Div u|^2\,dx+\dint
2\mu\left[\phi''(|u|)-\frac{\phi'(|u|)}{|u|}\right]|\nabla
|u||^2\,dx\\
&&\leq
-\dint\lambda\left[\phi''(|u|)-\frac{\phi'(|u|)}{|u|}\right]\sum_{ij}\partial_{i}u_j\frac{u_iu_j}{|u|^2}(\Div
u)\,dx+\dint \phi'(|u|)\frac{u}{|u|}\cdot (\rho F)\,dx\\
&&+\dint\phi'(|u|)\sum_{ij}\frac{u_j}{|u|}\partial_i(\rho
G_{ij})\,dx.
\end{eqnarray*}
We conclude noticing that, from the definition of $\nu$ in
(\ref{eq:h+Ng}), we have:
\begin{eqnarray*}
&&\dint \mu\frac{\phi'(|u|)}{|u|}|\nabla u|^2\,dx+\dint
\lambda\frac{\phi'(|u|)}{|u|}|\Div u|^2\,dx\\
&&\qquad \geq \dint \nu\frac{\phi'(|u|)}{|u|}|\nabla u|^2\,dx.
\end{eqnarray*}
\qed

In order to simplify the presentation, we will prove both theorems in
the framework of the second one. More precisely, we will show that
if  $\lambda=0$ then $u\in L^\infty_{\mathrm{loc}}([0,T]\times\R^N)$,
and that if $\lambda\neq0$ then we get some $L^p$ norms. This will
include Theorem \ref{theo:scalaire}, by applying the result
to $u=(\theta,\theta,\theta)$.

\section{Main propositions}
First, we introduce the function:
$$
v_k=[|u|-C_k]_+,
$$
where $C_k$ is an increasing sequence of positive numbers (to be
chosen later). Note that $v_k^2$ can be seen as a level set of
energy since $v_k^2=0$ for $|u|<C_k$ and is of order $|u|^2$ for
$|u|\gg C_k$.
We also consider a non-decreasing sequence of time $T_k$.

Then, we define
$$
U_k=\sup_{T_k<t<T}\left(\dint\rho(t,x)\frac{|v_k(t,x)|^2}{2}\,dx\right)
+\int_{T_k}^T\dint\nu|d_k(t,x)|^2\,dx\,dt,
$$
where:
$$
d_k^2=\frac{C_k\un_{\{|u|\geq
C_k\}}}{|u|}|\nabla|u||^2+\frac{v_k}{|u|}|\nabla u|^2.
$$

Note that if we take $C_k=0$, we get
$$
U_k=\sup_{T_k<t<T}\left(\dint\rho(t,x)\frac{|u(t,x)|^2}{2}\,dx\right)
+\int_{T_k}^T\dint\nu|\nabla u(t,x)|^2\,dx\,dt.
$$
Our goal is to show that for an appropriate choice of $C_k$,
the sequence $U_k$ goes to zero as $k$ goes to infinity. This will
be a consequence of the following propositions:

\begin{prop}[scalar case]\label{prop_NLsurcritique_scalaire}
 Consider $u$ solution to (\ref{eq_NSenergy}) with
$\lambda=0$ and let $F$ and $G$ verify
(\ref{eq:FGLinfty2})-(\ref{eq:pq}). For a fixed $K>0$, we define
$C_k$  by:
$$
C_k=K(1-2^{-k}).
$$
Then there exists $A_p\geq1$, $\beta_{1p}>1$, $\beta_{2p}>1$, and
$\gamma_p>0$, depending only on $p$ such that the following statements hold:
\begin{itemize}
\item If $u_0\in L^\infty(\R^3)$, then,
taking $T_k=0$ for all $k$ and  $K>2\|u_0\|_{L^\infty}$, we have:
$$
U_k\leq
\frac{A_p^k}{K^{\gamma_p}}\left(U_{k-1}^{\beta_{1p}}+U_{k-1}^{\beta_{2p}}\right)
\qquad \mbox{ for all } k\geq 2 .
$$
\item If $\rho\in L^\infty(0,T;L^r(\R^3))$ for some $r>3/2$, then for every $t_0>0$
we define $T_k=t_0(1-2^{-k})$, and there exists $A'_r\geq1$,
$\beta'_r>1$, and $\gamma_r'>0$ such that for any $K>0$:
$$
U_k \leq
\frac{A_p^k}{K^{\gamma_p}}U_{k-1}^{\beta_{1p}}+\frac{{A_r'}^k}{t_0K^{\gamma'_r}}
U_{k-1}^{\beta'_r}
\qquad \mbox{ for all }k\geq 2 .
$$
\end{itemize}

\end{prop}

\begin{prop}[system case]\label{prop_NLsurcritique_system}
Consider $u$ solution to (\ref{eq_NSenergy}) where $F$ and $G$
verify  (\ref{eq:FGLinfty2})-(\ref{eq:pq}), and assume that $\mu$
and $\lambda$ verify  (\ref{eq:h+Ng}) and (\ref{eq:lambda}). For a
fixed $K>0$, we define $C_k$ by:
$$
C_k=K2^k.
$$
Then there exists $0<\varepsilon<1$,  $\beta_{1p}>1$,
$\beta_{2p}>1$, and $\gamma_p>0$, depending only on $p$ such that
the following holds true:
\begin{itemize}
\item When $u_0\in L^\infty(\R^3)$, we take $T_k=0$ for all $k$. Then we have,
 for any $K>\|u_0\|_{L^\infty}$:
$$
U_k\leq
\frac{U_{k-1}^{\beta_{1p}}+U_{k-1}^{\beta_{2p}}}{K^{\gamma_p}}+\varepsilon\kappa
U_{k-1}  \qquad \mbox{ for all } k\geq2.
$$
\item When $\rho\in L^\infty(0,T;L^r(\R^3))$ for some  $r>3/2$, then for every $t_0>0$ we define $T_k=t_0(1-\eta^{-k})$.
Then there exists  $\beta'_r>1$, $\gamma_r'>0$ and
$0<\eta<1$ such that for any $K>0$  we have:
$$
U_k\leq
\frac{U_{k-1}^{\beta_p}}{K^{\gamma_p}}+\frac{U_{k-1}^{\beta'_r}}{(1-\eta)t_0K^{\gamma'_r}}
+\varepsilon\kappa U_{k-1} \qquad \mbox{ for all } k\geq2.
$$
\end{itemize}
\end{prop}


The next section is dedicated to the proofs of those two propositions.
In the following one we will show how those propositions indeed imply Theorem \ref{theo:scalaire} and Theorem \ref{theo:systeme}.

\section{Proof of Propositions \ref{prop_NLsurcritique_scalaire} and \ref{prop_NLsurcritique_system}}

This section is devoted to the proof of Propositions
\ref{prop_NLsurcritique_scalaire} and
\ref{prop_NLsurcritique_system}. The proof is split into several
steps.
\vspace{10pt}

\paragraph{Step 1: Evolution of $\dint \rho v_k^2\, dx$.}
The following lemma gives the inequality satisfied by the energy
of the level set function $v_k$:
\begin{lemm}\label{lemm_P}
Let $u$ be a solution of  (\ref{eq_NSenergy}) in
$Q=]0,T[\times\R^3$, then we have:
\begin{eqnarray}\nonumber
&&\qquad\frac{d}{dt}\dint \rho\frac{v_k^2}{2}\,dx+\dint \nu d_k^2\,dx\\
&&\leq -\dint\lambda r_k\,dx+\dint\frac{v_k}{|u|} u \cdot(\rho
F+\Div (\rho G))\,dx,\label{eq_NSvk}
\end{eqnarray}
where:
$$
r_k=
(\Div u)u\cdot\nabla|u|\frac{C_k}{|u|^2}\un_{\{|u|\geq C_k\}}.
$$
\end{lemm}
\pf
This lemma follows from (\ref{eq_NSenergy}) with
$$
\phi(y)=\frac{1}{2}(y-C_k)^2_+,
$$
and using the fact that
$$
\phi''(|u|)-\frac{\phi'(|u|)}{|u|}=\frac{C_k\unk}{|u|}.
$$
 \qed

\vspace{10pt}

\paragraph{ Step 2: First estimates on  $U_k$.}
Integrating  (\ref{eq_NSvk}) on $[\sigma,t]\times\Omega$, with
$T_{k-1}\leq \sigma\leq T_k\leq t\leq T$ we get
\begin{eqnarray*}
&&\dint\rho \frac{\left| v_k(t,x)\right|^2}{2}\,dx
+\int_{\sigma}^t\int  \nu d_k^2(s,x)\,dx\,ds\\
&&\qquad \leq \int \rho  \frac{|v_k(\sigma,x)|^2}{2}\,dx
-\int_\sigma^t\int \lambda r_{k} \,dx\,dt\\
&&\qquad\qquad+ \int_\sigma^t\int \frac{v_k}{|u|} u\cdot (\rho
F+\Div (\rho G) )\,dx\,dt.
\end{eqnarray*}
When $T_k=0$ and $C_{k}>\|u_0\|_{L^\infty}$ we have (with
$\sigma=T_k=0$):
\begin{eqnarray*}
&&\qquad\int \rho  \frac{|v_k(T_k,x)|^2}{2}\,dx\\
&&=\int \rho  \frac{|v_k(0,x)|^2}{2}\,dx\\
&&=\int \rho  \frac{(|u_0|-C_k)_+^2}{2}\,dx=0.
\end{eqnarray*}
When $T_k=t_0(1-\eta^{-k})$, then integrating with respect to $\sigma$ between
$T_{k-1}$ and $T_k$ and dividing by
$T_{k-1}-T_k=t_0\eta^{k-1}(1-\eta)$, we get
\begin{eqnarray*}
&&U_{k} = \sup_{t\in[T_k,T]}\left(\int\rho
\frac{\left|v_k(t,x)\right|^2}{2}\,dx
+\int_{T_k}^t\dint \nu d_k^2(s,x)\,dx\,ds\right)\\
&&\qquad \leq  \frac{1}{t_0\eta^{k-1}(1-\eta)}\int_{T_{k-1}}^{T_k}\int \rho \frac{|v_k(\sigma,x)|^2}{2}\,dx\,d\sigma\\
&&\qquad\qquad+\int_{T_{k-1}}^T\int \lambda(\rho) r_{k} \,dx\,dt\\
&&\qquad\qquad+ \int_{T_{k-1}}^T \int \frac{v_k}{|u|} u \cdot (\rho
F+\Div (\rho G))\,dx\,dt.
\end{eqnarray*}
The following lemma follows:
\begin{lemm}\label{lemm_eq_etakvk2}
\item[(i)] When $u_0\in L^\infty(\R^3)$,
we take $T_k=0$ for all $k$ and $C_{k}>\|u_0\|_{L^\infty}$, then we have:
\begin{equation}\label{eq_UkLinfty}
\begin{array}{l}
\displaystyle U_k\leq
 \int_{0}^T\int \lambda r_{k} \,dx\,dt\\
\displaystyle \qquad\qquad+ \int_{0}^T \int
\left(\frac{v_k}{|u|}\right)u\cdot (\rho F+\Div (\rho G))\,dx\,dt.
\end{array}
\end{equation}
\item[(ii)] When $u_0$ is only bounded in $L^2(\rho_0)$, we take
$T_k=t_0(1-\eta^{-k})$, then we have:
\begin{equation}\label{eq_UknoLinfty}
\begin{array}{l}
\displaystyle U_k\leq \frac{1}{t_0\eta^{k-1}(1-\eta)}\int_{T_{k-1}}^{T_k}\int \rho \frac{|v_k(\sigma,x)|^2}{2}\,dx\,d\sigma \\
\displaystyle \qquad\qquad+\int_{T_{k-1}}^T\int \lambda r_{k} \,dx\,dt\\
\displaystyle \qquad\qquad+ \int_{T_{k-1}}^T \int
\left(\frac{v_k}{|u|}\right)u\cdot (\rho F+\Div (\rho G))\,dx\,dt.
\end{array}
\end{equation}

\end{lemm}

\vspace{10pt}

\paragraph{ Step 3:  Some useful lemmas.}

We now want to show that the right-hand side in (\ref{eq_UkLinfty})
and (\ref{eq_UknoLinfty}) can be controled by terms of the form
$U_{k-1}^{\beta}$ with $\beta> 1$. This is the corner stone of De
Giorgi's method for the regularity of elliptic equation and a key step in this paper.

We start by giving the following technical lemma, which provides
some useful inequalities for the rest of the paper and  the proof of which
quite straightforward (and is given in the appendix for the comfort
of the reader):
\begin{lemm}\label{lemm_gradvk}
The function $u$ can be split in the following way:
$$
u=u\frac{v_k}{|u|}+u\left(1-\frac{v_k}{|u|}\right),
$$
where:
$$
\left|u\left(1-\frac{v_k}{|u|}\right)\right|\leq C_k.
$$
Moreover the following bounds hold:
\begin{eqnarray*}
&&\frac{v_k}{|u|}|\nabla u|\leq d_k,\\
&&\un_{\{|u|\geq C_{k}\}}|\nabla |u||\leq d_k,\\
&&|\nabla v_k|\leq d_k,\\
&&\left|\nabla\frac{uv_k}{|u|}\right|\leq 3d_k.
\end{eqnarray*}
\end{lemm}

The next Lemma will be crucial in what follows:
\begin{lemm}\label{lemm_tordu}
For every nonnegative  numbers $0<\alpha<1$, $\beta>1$ satisfying:
\begin{eqnarray*}
&&p_{1}=\frac{1}{\beta-\alpha}\geq1\\
&&q_{1}=\frac{3}{2\alpha+\beta}\geq1,
\end{eqnarray*}
we have:
$$
\|\rho^\alpha v_k^{2\beta}\|_{L^{p_{1}}(T_{k},T;L^{q_{1}}(\R^3))} \leq
U_{k}^{\beta}.
$$
\end{lemm}
\pf We obviously have:
$$
\|\rho^\alpha
v_k^{2\beta}\|_{L^{p_{1}}(T_{k},T;L^{q_{1}}(\R^3))}=\|\rho^{\alpha/\beta}
v_k^{2}\|^{\beta}_{L^{{p_{1}}\beta}(T_{k},T;L^{{q_{1}}\beta}(\R^3))}.
$$
Next, we note that:
$$
\rho^{\alpha/\beta} v_k^{2}=(\rho
v_k^{2})^{\alpha/\beta}v_k^{2(1-\frac{\alpha}{\beta})},
$$
and:
\begin{eqnarray*}
&& \|(\rho
v_k^{2})^{\alpha/\beta}\|_{L^\infty(T_k,T;L^{\beta/\alpha}(\R^3))}\leq
U_k^{\alpha/\beta},\\
&&\|v_k^{2(1-\frac{\alpha}{\beta})}\|_{L^{\frac{1}{1-\alpha/\beta}}
(T_k,T;L^{\frac{3}{1-\alpha/\beta}}(\R^3))}\leq \|\nabla
v_k\|^{2(1-\alpha/\beta)}_{L^2(]T_k,T[\times\R^3)}\leq
U_k^{1-\alpha/\beta}.
\end{eqnarray*}
(where we have used Lemma \ref{lemm_gradvk} for the last inequality).
It is now readily seen that H\"older's
inequalities give the result.\qquad\qed


\vspace{10pt}

\paragraph{ Step 4:  Terms involving $F$ and $G$.}

\begin{lemm}
If $F$ and $G$ satisfy conditions (\ref{eq:FGLinfty2})-(\ref{eq:pq}), then
there exists some $\beta>1$ such that
\begin{equation}\label{eq_GNL}
\left|-\int_{T_{k-1}}^T\dint  \frac{v_k}{|u|} u \cdot \Div (\rho
G)\,dx\,d\sigma\right|\leq C (C_k-C_{k-1})^{-\beta}
U_{k-1}^{(1+\beta)/2}
\end{equation}
and
\begin{equation}\label{eq_FNL}
\left|-\int_{T_{k-1}}^T\dint  \frac{v_k}{|u|} u \cdot \rho
F\,dx\,d\sigma\right|\leq C(C_k-C_{k-1})^{-2\beta+1} U_{k-1}^\beta.
\end{equation}
\end{lemm}

\pf First, we write
\begin{eqnarray*}
&&-\dint  \frac{v_k}{|u|} u \cdot \Div (\rho G)\,dx\\
&=&\int_{\R^3} \rho G:\left(\nabla u(1-\frac{C_k}{|u|})_++(u\otimes\nabla)|u|\frac{C_k}{|u|^2}\unk\right)\,dx\\
&\leq&C\left(\dint \rho^2\frac{|G|^2}{\nu}\unk\, dx\right)^{1/2}
\\
&&\qquad \times \left(\dint \nu(|\nabla u|^2(1-\frac{C_k}{|u|})_+^2+\nu
|\nabla|u||^2
\frac{C_k^2}{|u|^2}\unk)\, dx\right)^{1/2}\\
&\leq&C\left(\dint \frac{\rho^2|G|^2}{\nu}\unk\,
dx\right)^{1/2}\left(\dint \nu|d_k|^2\,dx\right)^{1/2}.
\end{eqnarray*}
Moreover, we have:
\begin{eqnarray*}
&& \int_{T_{k-1}}^{T} \into\frac{\rho^2|G|^2}{\nu}\unk\,dx\, d\sigma\\
&& \leq (C_k-C_{k-1})^{-2\beta}\int_{T_{k-1}}^{T}
\dint \rho^{1-\alpha} v_{k-1}^{2\beta}\frac{\rho^{1+\alpha}G^2}{\nu }\,dx\, d\sigma\\
&&\leq (C_k-C_{k-1})^{-2\beta}\left\|\frac{
\rho^{1+\alpha}|G|^2}{\nu}\right\|_{L^{p'_{1}}(T_{k-1},T;L^{q'_{1}}(\R^3))}
 \|\rho^{1-\alpha} v_k^{2\beta}\|_{L^{p_{1}}(T_{k-1},T;L^{q_{1}}(\R^3))} \\
&& \leq (C_k-C_{k-1})^{-2\beta}\left\|\frac{\rho^{1+\alpha}|G|^2}{\nu
}\right\|_{L^{p'_{1}}(T_{k-1},T;L^{q'_{1}}(\R^3))} U_{k-1}^\beta
\end{eqnarray*}
where we used Lemma \ref{lemm_tordu}. In order to have $\beta>1$, we need to take $p_{1}$ and $q_{1}$ such that
\begin{eqnarray*}
&&p_{1}< \frac{1}{\alpha}\\
&&q_{1}>\frac{3}{3-2\alpha}
\end{eqnarray*}
(using  Lemma \ref{lemm_tordu} with $1-\alpha$ instead of $\alpha$).
This leads to the following condition on $p'_{1}$ and $q'_{1}$ (conjugate of $p_{1}$ and $q_{1}$):
\begin{eqnarray*}
&&p'_{1}>\frac{1}{1-\alpha}\\
&&q'_{1}>\frac{3}{2\alpha}.
\end{eqnarray*}
It is readily seen that such coefficients will satisfy (\ref{eq:pq}), thus, using (\ref{eq:FGLinfty2}), we deduce
\begin{equation*}
\left|-\int_{T_{k-1}}^T\dint  \frac{v_k}{|u|} u \cdot \Div (\rho
G)\,dx\,d\sigma\right|\leq C (C_k-C_{k-1})^{-\beta}
U_{k-1}^{(1+\beta)/2}.
\end{equation*}
\vspace{5pt}

Next, we consider the term involving $F$. We have:
\begin{eqnarray*}
&&\left|-\int_{T_{k-1}}^T\dint  \frac{v_k}{|u|} u \cdot \rho F\,dx\,d\sigma\right|\\
&\leq&\int_{T_{k-1}}^T\dint  v_k \rho|F| \,dx\,d\sigma\\
&\leq&(C_k-C_{k-1})^{-(2\beta-1)}\int_{T_{k-1}}^{T} \dint
\rho^{1-\alpha} v_{k-1}^{2\beta}\frac{\rho^\alpha
F}{\rho^\alpha}\,dx\,
d\sigma\\
&\leq&(C_k-C_{k-1})^{-2\beta+1}\left\|\rho^\alpha|F|
\right\|_{L^{p'_{1}}(T_{k-1},T;L^{q'_{1}}(\R^3))}
 \|\rho^{1-\alpha} v_k^{2\beta}\|_{L^{p_{1}}(T_{k-1},T;L^{q_{1}}(\R^3)} \\
&\leq&(C_k-C_{k-1})^{-2\beta+1}\left\|\rho^\alpha|F|
\right\|_{L^{p'_{1}}(T_{k-1},T;L^{q'_{1}}(\R^3))} U_{k-1}^\beta
\end{eqnarray*}
where we used lemma \ref{lemm_tordu} as before (with the same conditions on $p_{1}$ and $q_{1}$).
Using (\ref{eq:FGLinfty2}),
we deduce
\begin{equation*}
\left|-\int_{T_{k-1}}^T\dint  \frac{v_k}{|u|} u \cdot \rho
F\,dx\,d\sigma\right|\leq C(C_k-C_{k-1})^{-2\beta+1} U_{k-1}^\beta.
\end{equation*}

\noindent{\bf Remark:} Note that when  $\lambda = 0$ and $u_{0}\in L^{\infty}(\R^{3})$, and if
we take $K>2\|u_0\|_{L^\infty}$, $T_k=0$ and $C_k=K(1-2^{-k})$, then
(\ref{eq_UkLinfty}), (\ref{eq_GNL}) and (\ref{eq_FNL})  yield
$$ U_{k} \leq C (C_{k}-C_{k-1})^{-\gamma'} U_{k-1}^{\beta'}$$
with $\gamma'=\sup(\beta,2\beta-1)$ and $\beta'=\inf(\beta,(1+\beta)/2)>1$.
Since $C_k-C_{k-1}=K 2^{-k}$, we deduce
$$ U_{k} \leq C K^{-\gamma'} 2^{\gamma' k} U_{k-1}^{\beta'}.$$
This gives the first part of Proposition \ref{prop_NLsurcritique_scalaire}.
In the next step we will show how to deal with unbounded initial data.

\vspace{10pt}

\paragraph{ Step 5: Control of the time layer (case
$T_k\neq 0$).}

\begin{lemm}\label{lemm_layer_t}
If $\rho$ is bounded in $L^{\infty}(0,T;L^{{r} }(\R^3))$ for some
${r} >3/2$, then we have:
\begin{eqnarray*}
&&\qquad\frac{1}{t_0\eta^{k-1}(1-\eta)}\int_{T_{k-1}}^{T_k}\int \rho
\frac{|v_k(\sigma,x)|^2}{2}\,dx\,d\sigma \\
&&\leq
(t_0\eta^{k-1}(1-\eta))^{-1}(C_{k}-C_{k-1})^{-\frac{2\alpha}{3}}
||\rho||^{\frac{3-\alpha}{3}}_{L^{\infty}(0,T;L^{{r} }(\R^3))}
U^{1+\alpha/3}_{k-1}
\end{eqnarray*}
 with
$$ \alpha = \frac{2 {r} -3}{{r} -1}>0$$
\end{lemm}
{\it Proof.} First, using H\"older and Sobolev inequalities, we
have:
\begin{eqnarray*}
 && \int_{T_{k-1}}^{T_k}\int \rho \frac{|v_k(\sigma,x)|^2}{2}\,dx\,d\sigma \\
 && \leq\left[ \int_{T_{k-1}}^{T_k}\left( \int |v_k(\sigma,x)|^6 \,dx\right)^{1/3} \,d\sigma\right] \sup_{t \in[T_{k-1},T_{k}]} \left( \int \rho ^{3/2}(x,t) \mathrm{1}_{\{v_{k}>0\}}\,dx\right)^{2/3} \\
 && \leq\left[ \int_{T_{k-1}}^{T} \int |\na v_k(\sigma,x)|^2 \,dx \,d\sigma\right]  \sup_{t \in[T_{k-1},T_{k}]} \left( \int \rho ^{3/2}(x,t) \mathrm{1}_{\{v_{k}>0\}}\,dx \right)^{2/3}\\
  && \leq C U_{k-1} \, \sup_{t \in[T_{k-1},T_{k}]} \left( \int \rho ^{3/2}(x,t) \mathrm{1}_{\{v_{k}>0\}}\,dx
  \right)^{2/3}.
 \end{eqnarray*}
Next, we note that for all $x$ such that  $v_k(x)>0$, we have
$|u(x)|>C_{k}$ and so
\begin{eqnarray*}
v_{k-1}(x)&=&[|u(x)|-C_{k-1}]_+\\
&=&[|u(x)|-C_{k}+(C_{k}-C_{k-1})]_+\\
&>&C_{k}-C_{k-1},
\end{eqnarray*}
which yields
$$  \mathrm{1}_{\{v_{k}>0\}}  \leq (C_{k}-C_{k-1})^{-1}v_{k-1} .$$
It follows that
\begin{eqnarray*}
&& \left( \int \rho ^{\frac{3}{2}}(x,t) \mathrm{1}_{\{v_{k}>0\}}\,dx \right)^{\frac{2}{3}} \\
&& \leq  (C_{k}-C_{k-1})^{-\frac{2\alpha}{3}} \left( \int \rho ^{\frac{3-\alpha}{2}}(x,t) \rho^{\frac{\alpha}{2}}(x,t) |v_{k-1}|^{\alpha}\,dx \right)^{\frac{2}{3}} \\
& & \leq  (C_{k}-C_{k-1})^{-\frac{2\alpha}{3}} \left( \int \rho
^{\frac{3-\alpha}{2-\alpha}}(x,t)\,
dx\right)^{\frac{2-\alpha}{3}}\left(\int  \rho(x,t)
|v_{k-1}|^{2}\,dx \right)^{\frac{\alpha}{3}}
\end{eqnarray*}
With $\alpha = \frac{2{r} -3}{{r} -1}$ (so that
$\frac{3-\alpha}{2-\alpha} = {r} $), we deduce
\begin{eqnarray*}
&&\sup_{t\in[T_{k-1},T_{k}]} \left( \int \rho ^{\frac{3}{2}}(x,t) \mathrm{1}_{\{v_{k}>0\}}\,dx \right)^{\frac{2}{3}} \\
& & \leq  (C_{k}-C_{k-1})^{-2\alpha/3} ||\rho
||^{(3-\alpha)/3}_{L^{\infty}(0,T;L^{{r} }(\R^3))}
 \left(\sup_{t\in[T_{k-1},T_{k}]}  \int  \rho(x,t) |v_{k-1}|^{2}\,dx \right)^{\frac{\alpha}{3}}
\\
& & \leq  (C_{k}-C_{k-1})^{-2\alpha/3} ||\rho
||^{\frac{{r} }{3({r} -1)}}_{L^{\infty}(0,T;L^{{r} }(\R^3))}
U_{k-1}^{\frac{2{r} -3}{3({r} -1)}}
\end{eqnarray*}
\qed

\noindent{\bf Remark:} When $\lambda=0$, we can take $T_k=1-2^{-k}$
and  $C_k=K(1-2^{-k})$. Proceeding as before, (\ref{eq_UknoLinfty}), (\ref{eq_GNL}), (\ref{eq_FNL}) and Lemma \ref{lemm_layer_t} give the second part of Proposition \ref{prop_NLsurcritique_scalaire}.

\vspace{10pt}

\paragraph{ Step 6: The second viscosity term.}
It  only remains to  control the term corresponding to the second viscosity coefficient. This is achieved by the following lemma:
\begin{lemm}
Under the assumptions (\ref{eq:h+Ng}) and (\ref{eq:lambda}), we have
$$
\left|\int_{T_{k-1}}^T\int_{\R^3}\lambda
r_k\,dx\,dt\right|\\
\leq\kappa \left(\frac{1}{3}
\frac{C_k}{C_k-C_{k-1}}\right)^{1/2}U_{k-1}.
$$
\end{lemm}
\pf We have
\begin{eqnarray*}
&&\qquad \left|\int_{T_{k-1}}^T\int_{\R^3}\lambda
r_k\,dx\,dt\right|\\
&&= \left|\int_{T_{k-1}}^T\int_{\R^3}\lambda
(\Div u)u\cdot\nabla|u|\frac{C_k}{|u|^2}\un_{\{|u|\geq C_k\}}
\,dx\,dt\right|\\
&&\leq
\left(\int_{T_{k-1}}^T\int_{\R^3}|\lambda|\frac{v_{k-1}}{|u|}|\Div
u|^2\,dx\,dt\right)^{1/2}\left(\int_{T_{k-1}}^T\int_{\R^3}
|\lambda|\frac{C_k^2}{|u|v_{k-1}}|\nabla|u||^2\unk\right)^{1/2}\\
&&\leq \left(\kappa\int_{T_{k-1}}^T\int_{\R^3}\nu d_k^2
\,dx\,dt\right)^{1/2}\left(\int_{T_{k-1}}^T\int_{\R^3}
\frac{|\lambda|}{\nu}\frac{\nu
C_k}{v_{k-1}}\unk|d_k|^2\,dx\,dt\right)^{1/2}.
\end{eqnarray*}
Noticing that
$$
\frac{C_k}{v_{k-1}}\unk\frac{|\lambda|}{\nu}\leq \frac{\kappa}{3}
\frac{C_k}{C_k-C_{k-1}},
$$
we deduce
$$
\left|\int_{T_{k-1}}^T\int_{\R^3}\lambda
r_k\,dx\,dt\right|\\
\leq\kappa \left(\frac{1}{3}
\frac{C_k}{C_k-C_{k-1}}\right)^{1/2}U_{k-1}.
$$

\vspace{10pt}


\vspace{10pt}

\paragraph{ Step 7: Conclusion.}

We have already proven Proposition
\ref{prop_NLsurcritique_scalaire}, which corresponds to $\lambda=0$ (the first part was completed at
the end of Step 4, and the second part at the end of Step 5).
In order to control the term due to the second viscosity coefficient, and prove Proposition \ref{prop_NLsurcritique_system}, we need the quantity $\frac{C_k}{C_k-C_{k-1}}$ to be bounded.
We thus take  $C_k=K2^k$.
Then we have
$$
\frac{C_k}{C_k-C_{k-1}}\leq2,
$$
and so:
\begin{equation}
\left|\int_{T_{k-1}}^T\int_{\R^3}\lambda r_k\,dx\,dt\right|\\
\leq \varepsilon\kappa U_{k-1},\label{eq_r_k}
\end{equation}
where $\varepsilon=\sqrt{2/3}<1$.
The first part of Proposition \ref{prop_NLsurcritique_system} now follows using
(\ref{eq_UkLinfty}), (\ref{eq_FNL}),  (\ref{eq_GNL}), and
(\ref{eq_r_k}).

The second part follows (\ref{eq_UknoLinfty}),
(\ref{eq_FNL}), (\ref{eq_GNL}), Lemma \ref{lemm_layer_t},
(\ref{eq_r_k}) and the fact that since
$$
C_k-C_{k-1}=K2^{k-1}>1,
$$
we can always choose $\eta<1$ such that
$$
\eta^{-k+1}(C_k-C_{k-1})^{-\frac{2\alpha}{3}}\leq 1 \qquad \mbox{ for all $k \geq 2$}.
$$

\vspace{20pt}

\section{Proofs of the theorems}
We can now complete the proofs of the theorems.
We have to show that the sequences $U_{k}$ constructed in the previous section converge to zero as $k$ goes to infinity for an appropriate choice of constant $K$.
It will be a
consequence of the following lemmas:
\begin{lemm}\label{lemm_technique_scalar}
Let ${U_k}$ be a sequence satisfying
\begin{eqnarray*}
&&0\leq U_0\leq C,\\
&&0\leq U_{k}\leq
\frac{A^k}{K}(U_{k-1}^{\beta_1}+U_{k-1}^{\beta_2}),\qquad \forall
k\geq1,
\end{eqnarray*}
for some constants $A\geq1$, $1<\beta_1<\beta_2$ and $C>0$.

Then there exists $K_0$ such that for every $K>K_0$ the
 sequence $U_k$ converges to 0 when $k$ goes to
infinity.
\end{lemm}

\begin{lemm}\label{lemm_technique_system}
Let  ${U_k}$ be a sequence satisfying:
\begin{eqnarray*}
&&0\leq U_0\leq C,\\
&&0\leq U_{k}\leq
\frac{1}{K}(U_{k-1}^{\beta_1}+U_{k-1}^{\beta_2})+\varepsilon\kappa
U_{k-1},\qquad \forall k\geq1
\end{eqnarray*}
for some constants $0<\varepsilon<1$, $0<\kappa<1$,
$1<\beta_1<\beta_2$ and $C>0$.
Then there exists $K_0$ such that for every $K>K_0$ the
 sequence $U_k$ converges to 0 when $k$ goes to
infinity. Moreover there exists $0<\varepsilon_1<1$  and $C_K>0$
such that:
$$
U_{k}\leq C_K (\varepsilon_1\kappa)^k.
$$
\end{lemm}
\noindent{\bf Proof of Lemma \ref{lemm_technique_scalar}:} We introduce
$$
\overline{U}_k=\frac{U_k}{K^{\frac{1}{\beta_1-1}}}.
$$
As long as $K>1$ we get for every $k$:
\begin{eqnarray*}
\overline{U}_k&\leq&
A^k\left(\overline{U}^{\beta_1}_{k-1}+\frac{\overline{U}^{\beta_2}_{k-1}}
{K^{\frac{\beta_2-\beta_1}{\beta_1-1}}}\right)\\
&\leq&A^k\left(\overline{U}^{\beta_1}_{k-1}+\overline{U}^{\beta_2}_{k-1}
\right).
\end{eqnarray*}
Next, we consider the sequence
$$
W_k=(2A)^kW_{k-1}^{\beta_1},
$$
and we claim that if $W_0$ is small enough then $W_k<1$ for every $k$
and $W_{k}$ converges to 0. Indeed, introducing
$$
\overline{W}_k=(2A)^{\frac{k}{\beta_1-1}}(2A)^{\frac{1}{(\beta_1-1)^2}}W_k,
$$
we get
$$
 0\leq\overline{W}_{k+1}\leq\overline{W}_{k}^{\beta_1}.
$$
So if $W_0\leq C^*_0=(2A)^{-1/(\beta-1)^2}$, we have
$\overline{W}_0\leq1$ and by induction $\overline{W}_k\leq 1$ for
every $k$. This gives:
$$
W_k\leq (2A)^{\frac{-k}{\beta-1}}C^{\frac{-1}{(\beta-1)^2}}.
$$
Since $2A>1$,
 this shows
that $W_k$ converges to 0 when k goes to infinity. Therefore,  for
$K$ big enough, $\overline{U}_0\leq W_0$ and as long as
$\overline{U}_{k-1}<1$, we have $\overline{U}_k\leq W_k$. So by
induction we show that this is valid for every $k$ and so
$\overline{U}_k$ converges to 0. This implies that $U_k$ converges
to 0.\qquad\qed

\vskip1cm \noindent{\bf Proof  of Lemma
\ref{lemm_technique_system}:} Take $\varepsilon_1$ such that:
$$\varepsilon<\varepsilon_1<1 $$ and $K$ big enough such that
$$
\frac{C^{\beta_1-1}+C^{\beta_2-1}}{K}\leq
(\varepsilon_1-\varepsilon)\kappa.
$$
Then $U_1\leq \varepsilon_1 \kappa U_0$ and we can  show by
induction that
$$
U_{k}\leq(\varepsilon_1\kappa)^{k}C.
$$
\qed
\vspace{10pt}

Theorem \ref{theo:scalaire} (respectively Theorem \ref{theo:systeme})  is now a straighforward consequence of Proposition \ref{prop_NLsurcritique_scalaire} (resp. Proposition \ref{prop_NLsurcritique_system}) and Lemma \ref{lemm_technique_scalar} (resp. Lemma \ref{lemm_technique_system}).
 \vspace{5pt}

\noindent{\bf Proof of Theorem \ref{theo:scalaire}:}
Let $u=(\theta,\theta,\theta)$. Then $u$ verifies (\ref{eq_NSenergy}) with $\lambda=0$, so  we can use Proposition \ref{prop_NLsurcritique_scalaire}.
We can then use Lemma
\ref{lemm_technique_scalar} with $K= \inf
(K^{\gamma_p},K^{\gamma'_p})$ and $A=\sup
(A_p,A'_q)$.
If follows that for $K$ big enough $U_k$
converges to 0.

Next, we observe that
\begin{eqnarray*}
U_k&\geq&\int_0^\infty\dint\un_{\{t\in (T_k,T)\}}d_k^2\,dx\,dt.
\end{eqnarray*}
So when $T_k=0$, we get that  $d_k^2\leq 2|\nabla u|^2$ and $d_k$ converges
almost everywhere to $d_\infty$ defined with $C_\infty=K$.
By Lebesgue's dominated convergence theorem  and Lemma
\ref{lemm_gradvk} we deduce
$$
\int_0^T\dint |\nabla(|u|-K)_+|^2\,dx\,dt\leq\int_0^T\dint
d_\infty^2\,dx\,dt=0,
$$
so $(|u|-K)_+=0$ and $|\theta|\leq K/3$. If $T_k=t_0(1-2^{-k})$, we proceed in
the same way and we find that:
$$
\int_{t_0}^T\dint |\nabla(|u|-K)_+|^2\,dx\,dt\leq\int_0^T\dint
d_\infty^2\,dx\,dt=0.
$$
It follows that  $\theta$ is bounded on $(t_0,T)\times\R^3$, and this is true for
any $t_0>0$.\qquad\qed

\vskip1cm \noindent{\bf Proof of Theorem \ref{theo:systeme}:} We now use
Proposition \ref{prop_NLsurcritique_system}
and Lemma
\ref{lemm_technique_system} with $K=\inf
(K^{\gamma_p},K^{\gamma'_p})$ .
It follows that for $K$ big enough we get
$$
U_k\leq C_K (\varepsilon_1\kappa)^k.
$$
>From Lemma \ref{lemm_gradvk}, we have:
$$
U_k\geq \|v_k\|^2_{L^2(\sup T_k,T;L^6(\R^3))}\geq
\|v_k\|^2_{L^2(\sup T_k,T;L^2_{\mathrm{loc}}(\R^3))},
$$
so using Tchebichev's inequality, we find that for every $R>0$ we have
\begin{eqnarray*}
{\mathcal{L}}(\{|u|\geq2K 2^k\}\cap\{(\sup T_k,T)\times
B(R)\})&\leq& {\mathcal{L}}(\{|u|\geq2C_k\}\cap\{(\sup T_k,T)\times
B(R)\})\\
&\leq&\frac{U_k}{C_k^2}\\
&\leq&C_K 2^{-2k}(\varepsilon_1\kappa)^k\\
&\leq&C_K 2^{-2k-pk}
\end{eqnarray*}
for some $p>\log_2\frac{1}{\kappa}$.
This implies that $u$ lies in
$L^{2+p,*}_{\mathrm{loc}}$, and thus is bounded in
$L^{2+\log_2\frac{1}{\kappa}}_{\mathrm{loc}}$.\qquad\qed

\appendix
\section{Proof of Lemma \ref{lemm_gradvk}}
\pf The function $(1-v_k/|u|)$ is Lipshitz and equal to:
\begin{eqnarray*}
1-\frac{v_k}{|u|}&=&1\qquad\qquad\mathrm{if} \qquad
|u|\leq C_k\\
&=&\frac{ C_k}{|u|}\qquad\qquad\mathrm{if} \qquad |u|\geq C_k.
\end{eqnarray*}
Therefore:
$$
\left|u\left(1-\frac{v_k}{|u|}\right)\right|\leq C_k.
$$
Let us first show  that:
\begin{eqnarray}\label{eq_nablaudk}
&&\frac{v_k}{|u|}|\nabla u|\leq d_k\\
\label{eq_nabla|u|dk} &&\un_{\{|u|\geq C_k\}}|\nabla |u||\leq d_k.
\end{eqnarray}
Statement (\ref{eq_nablaudk}) comes from the definition of $d_k$ and
the fact that $v_k\leq |u|$:
$$
d_k^2\geq \frac{v_k}{|u|}|\nabla u|^2\geq
\left(\frac{v_k}{|u|}|\nabla u|\right)^2.
$$
To show (\ref{eq_nabla|u|dk}), notice that:
$$
|\nabla|u||^2=\left|\frac{u}{|u|}\nabla u\right|^2\leq |\nabla u|^2.
$$
So:
$$
d_k^2\geq\frac{( C_k)\un_{\{|u|\geq (
C_k)\}}+v_k}{|u|}|\nabla|u||^2,
$$
with:
$$
(( C_k)+v_k)\un_{\{|u|\geq ( C_k)\}}=|u|\un_{\{|u|\geq ( C_k)\}}.
$$
So:
$$
d_k^2\geq\un_{\{|u|\geq ( C_k)\}}|\nabla|u||^2.
$$
Then the bound on $\nabla v_k$ follows (\ref{eq_nabla|u|dk}) since:
$$
|\nabla v_k|=|\nabla|u||\un_{\{|u|\geq ( C_k)\}}.
$$
To find the last inequality we fist write:
$$
\nabla \left( \frac{uv_k}{|u|}\right)=\frac{u}{|u|}\nabla v_k
+v_k\nabla\left(\frac{u}{|u|}\right).
$$
The first term can be bounded by:
$$
\left|\frac{u}{|u|}\nabla v_k\right|\leq |\nabla v_k|\leq d_k.
$$
The second one can be rewritten in the following way:
$$
v_k\nabla\left(\frac{u}{|u|}\right)=\frac{v_k}{|u|}\nabla
u-\frac{v_k u}{|u|^2}\nabla|u|.
$$
So, thanks to (\ref{eq_nablaudk}) and (\ref{eq_nabla|u|dk}):
\begin{eqnarray*}
\left|v_k\nabla\left(\frac{u}{|u|}\right)\right|&\leq&\frac{v_k}{|u|}|\nabla
u|+ \un_{\{|u|\geq ( C_k)\}}|\nabla |u||\\
&\leq&2 d_k.
\end{eqnarray*}
This gives:
$$
\left|\nabla \left( \frac{uv_k}{|u|}\right)\right|\leq 3d_k.
$$
This ends the proof of the lemma.\qquad\qed

\bibliography{biblio}

\begin{thebibliography}{10}

\bibitem{BDTemp}
Didier Bresch and Benoit Desjardins.
\newblock On the existence of global weak solutions to the navier-stokes
  equations for viscous compressible and heat conducting fluids.
\newblock {\em Preprint}.

\bibitem{BD2}
Didier Bresch and Beno{\^{\i}}t Desjardins.
\newblock Existence of global weak solutions for a 2{D} viscous shallow water
  equations and convergence to the quasi-geostrophic model.
\newblock {\em Comm. Math. Phys.}, 238(1-2):211--223, 2003.

\bibitem{BD3}
Didier Bresch and Benoit Desjardins.
\newblock Some diffusive capillary models of korteweg type.
\newblock {\em Notes aux comptes-rendus de l'academie des Sciences,
  Mechanique}, 331, 2003.

\bibitem{CKN}
L.~Caffarelli, R.~Kohn, and L.~Nirenberg.
\newblock Partial regularity of suitable weak solutions of the
  {N}avier-{S}tokes equations.
\newblock {\em Comm. Pure Appl. Math.}, 35(6):771--831, 1982.

\bibitem{DeGiorgi}
Ennio De~Giorgi.
\newblock Sulla differenziabilit\`a e l'analiticit\`a delle estremali degli
  integrali multipli regolari.
\newblock {\em Mem. Accad. Sci. Torino. Cl. Sci. Fis. Mat. Nat. (3)}, 3:25--43,
  1957.

\bibitem{F3}
Eduard Feireisl.
\newblock On the motion of a viscous, compressible, and heat conducting fluid.
\newblock {\em Indiana Univ. Math. J.}, 53(6):1705--1738, 2004.

\bibitem{F2}
Eduard Feireisl, Anton{\'{\i}}n Novotn{\'y}, and Hana Petzeltov{\'a}.
\newblock On the existence of globally defined weak solutions to the
  {N}avier-{S}tokes equations.
\newblock {\em J. Math. Fluid Mech.}, 3(4):358--392, 2001.

\bibitem{Hoff}
David Hoff.
\newblock Global existence for {$1$}{D}, compressible, isentropic
  {N}avier-{S}tokes equations with large initial data.
\newblock {\em Trans. Amer. Math. Soc.}, 303(1):169--181, 1987.

\bibitem{Hoff3}
David Hoff.
\newblock Global solutions of the {N}avier-{S}tokes equations for
  multidimensional compressible flow with discontinuous initial data.
\newblock {\em J. Differential Equations}, 120(1):215--254, 1995.

\bibitem{Hoff2}
David Hoff.
\newblock Strong convergence to global solutions for multidimensional flows of
  compressible, viscous fluids with polytropic equations of state and
  discontinuous initial data.
\newblock {\em Arch. Rational Mech. Anal.}, 132(1):1--14, 1995.

\bibitem{J1}
Song Jiang and Ping Zhang.
\newblock Axisymmetric solutions of the 3{D} {N}avier-{S}tokes equations for
  compressible isentropic fluids.
\newblock {\em J. Math. Pures Appl. (9)}, 82(8):949--973, 2003.

\bibitem{K}
A.~V. Kazhikhov and V.~V. Shelukhin.
\newblock Unique global solution with respect to time of initial-boundary value
  problems for one-dimensional equations of a viscous gas.
\newblock {\em Prikl. Mat. Meh.}, 41(2):282--291, 1977.

\bibitem{PLL2}
Pierre-Louis Lions.
\newblock {\em Mathematical topics in fluid mechanics. {V}ol. 2}, volume~10 of
  {\em Oxford Lecture Series in Mathematics and its Applications}.
\newblock The Clarendon Press Oxford University Press, New York, 1998.
\newblock Compressible models, Oxford Science Publications.

\bibitem{nishida}
Akitaka Matsumura and Takaaki Nishida.
\newblock The initial value problem for the equations of motion of compressible
  viscous and heat-conductive fluids.
\newblock {\em Proc. Japan Acad. Ser. A Math. Sci.}, 55(9):337--342, 1979.

\bibitem{MV}
Antoine Mellet and Alexis Vasseur.
\newblock On the compressible barotropic {N}avier-{S}tokes equations.
\newblock {\em Preprint}, 2005.

\bibitem{Serre}
Denis Serre.
\newblock Solutions faibles globales des \'equations de {N}avier-{S}tokes pour
  un fluide compressible.
\newblock {\em C. R. Acad. Sci. Paris S\'er. I Math.}, 303(13):639--642, 1986.

\bibitem{Vaigant}
V.~A. Va{\u\i}gant.
\newblock An example of the nonexistence with respect to time of the global
  solution of {N}avier-{S}tokes equations for a compressible viscous barotropic
  fluid.
\newblock {\em Dokl. Akad. Nauk}, 339(2):155--156, 1994.

\bibitem{Vasseur}
Alexis Vasseur.
\newblock A new proof of partial regularity of solutions to {N}avier-{S}tokes
  equations.
\newblock {\em To appear in NoDEA}.

\end{thebibliography}

\end{document}